\documentclass[11pt,a4paper,leqno]{amsart}

\usepackage[margin=1.2in]{geometry}
\usepackage{amsmath}
\usepackage{latexsym}
\usepackage{amsfonts}
\usepackage{amssymb}
\usepackage{amsmath}
\usepackage{amsthm}
\usepackage{amsxtra}
\usepackage{amscd}
\usepackage[foot]{amsaddr}
\usepackage{a4wide}
\usepackage{color}
\usepackage[utf8]{inputenc}
\usepackage{todonotes}
\usepackage{mathrsfs}
\usepackage{hyperref}
\usepackage{cite}
\usepackage{enumerate}

\usepackage{url}
\usepackage{orcidlink}
\usepackage{csquotes}
\numberwithin{equation}{section}

\usepackage{etoolbox}

\newcommand{\Sch}{\mathscr{S}}
\newcommand{\R}{\mathbb{R}}
\newcommand{\N}{\mathbb{N}}

\makeatletter
\patchcmd{\@settitle}
  {\uppercasenonmath\@title}
  {}{}{}

\patchcmd{\@setauthors}
  {\MakeUppercase{\authors}}
  {\authors}{}{}
\makeatother

\title[Erratum \enquote{Power boundedness for weighted composition operators}]{Erratum to \enquote{Power boundedness and related properties for weighted composition operators on $\Sch(\R^{d})$}\\ Journal of Functional Analysis 288 (2025), 110745}

\author[V.~Asensio]{Vicente Asensio$^1$ \orcidlink{0000-0002-7615-4140}}

\address{$^1$Instituto Universitario de Matem\'atica Pura y Aplicada IUMPA, Universitat Polit\`ecnica de Val\`encia, Camino de Vera, s/n, E-46022 Valencia, Spain}
\email{viaslo@upv.es}

\author[E.~Jord\'a]{Enrique Jord\'a$^2$ \orcidlink{0000-0003-2980-1699}}

\address{$^2$Departamento de Matem\'atica Aplicada, E. Polit\'{e}cnica Superior
	de Alcoy, Universidad Polit\'ecnica de Valencia, Plaza Ferr\'andiz
	y Carbonell 2, E-03801 Alcoy (Alicante), Spain}
\email{ejorda@mat.upv.es}
 
\author[T.~Kalmes]{Thomas Kalmes$^3$ \orcidlink{0000-0001-7542-1334}}

\address{$^3$Faculty of Mathematics, Chemnitz University of Technology, 09107 Chemnitz, Germany}
\email{thomas.kalmes@math.tu-chemnitz.de}

\date{}
\begin{document}

\begin{abstract}
As was kindly pointed out to us by S.~Li, the proofs of Theorem~2.11 and Lemma~3.1 in our original paper contain errors that require modifications of the corresponding statements. We also correct a gap in the proof of Lemma~2.1, whose statement remains valid. Furthermore, we explain the consequences for several subsequent results and provide the necessary corrections to their statements and proofs.
\end{abstract}

\maketitle

\section*{Introduction}

In \cite{AJK}, we studied power boundedness and related properties of weighted composition operators $C_{\psi,\varphi}$ on the Schwartz space $\mathcal{S}(\mathbb{R}^d)$, where $C_{\psi,\varphi}(f)=\psi\cdot(f\circ\varphi)$ for $f\in\mathscr{S}(\mathbb{R}^d)$, with $\psi\in C^\infty(\R^d)$ and a smooth mapping $\varphi\colon\R^d\rightarrow\R^d$, the symbol of $C_{\psi,\varphi}$.

As kindly pointed out to us by S.~Li, the proofs of Theorem~2.11 and Lemma~3.1 in \cite{AJK} contain errors. We are grateful to S.~Li for bringing these issues to our attention and for further investigating their consequences \cite{LiCorrections}.

In the present erratum, we correct Theorem~2.11 by imposing an additional boundedness assumption on the mapping $\varphi$, and Lemma~3.1 by strengthening the original condition (3.3). We also provide a new proof of Lemma~2.1, whose statement remains valid despite a gap in the original proof.

These corrections affect several subsequent results. In particular, we amend Theorems~4.5 and~4.12, whose corrected statements correspond to Remarks~4.6 and~4.13, respectively, of \cite{AJK}. Furthermore, we remove an incorrect assertion from Proposition~5.4, withdraw Proposition~5.5, and modify Theorems 5.8 and 5.10 accordingly. We also provide new proofs where required.




We are grateful to S.~Li for bringing the errors in the proofs of Theorem~2.11 and Lemma~3.1 to our attention, and for sharing his further analysis with us \cite{LiCorrections}. In this erratum, we restrict ourselves to correcting the affected statements and proofs, without reproducing his additional results.

\section{Corrected proof of Lemma 2.1}
The proof of Lemma~2.1 in \cite{AJK} incorrectly assumes that the points at which the relevant coefficient becomes large can be chosen to tend to infinity. This does not follow from hypothesis~\cite[(2.3)]{AJK}. The lemma itself remains valid, as the following proof shows.

\begin{proof}[Replacement proof of Lemma 2.1]
Fix $\alpha,\lambda\in\N_0^d$ with $|\lambda|\le |\alpha|$ and $p>0$ satisfying~\cite[(2.3)]{AJK}. Suppose, by contradiction, that for every $f\in\Sch(\R^d)$,
\[
\sup_{i\in I}\sup_{x\in\R^d}(1+|x|)^p \bigl| (\psi_i \cdot (f \circ \varphi_i))^{(\alpha)}(x) \bigr|<\infty.
\]
For $i\in I$ and $x\in\R^d$ define the continuous linear functional
\[
\Lambda_{i,x}(f)=(1+|x|)^p \big(\psi_i \cdot (f \circ \varphi_i)\big)^{(\alpha)}(x),\qquad f\in\Sch(\R^d).
\]
By assumption, the family $(\Lambda_{i,x})_{i,x}$ is pointwise bounded. Since $\Sch(\R^d)$ is a Fr\'echet space, the uniform boundedness principle gives $C>0$ and an integer $m\ge |\alpha|$ such that
\begin{equation}\label{eq:UBP}
|\Lambda_{i,x}(f)|\le C\|f\|_m\qquad(i\in I,\ x\in\R^d,\ f\in\Sch(\R^d)),
\end{equation}
where
\[
\|f\|_m=\max_{|\beta|\le m}\sup_{y\in\R^d}(1+|y|)^m|f^{(\beta)}(y)|.
\]
Choose $g\in \mathscr{D}(B(0,1))$ such that
\[
g^{(\lambda)}(0)=1,\qquad
g^{(\beta)}(0)=0\quad (|\beta|\le |\alpha|,\ \beta\ne\lambda).
\]
For each $i \in I, x \in \mathbb{R}^d$ put $f_{i,x}(y)=g(y-\varphi_i(x))$. Translation yields a constant $C_m$ independent of $i,x$ for which
\begin{equation}\label{eq:translate}
\|f_{i,x}\|_m\le C_m(1+|\varphi_i(x)|)^m.
\end{equation}
By formula \cite[(2.2)]{AJK} and the prescribed jet of $g$,
\[
\bigl(\psi_i \cdot (f_{i,x}\circ\varphi_i)\bigr)^{(\alpha)}(x)
=F_{\alpha,\lambda}^{\varphi_i,\psi_i}(x).
\]
Applying \eqref{eq:UBP} to $f_{i,x}$ and then \eqref{eq:translate}, we obtain
\[
\sup_{i\in I}\sup_{x\in\R^d}
\frac{(1+|x|)^p}{(1+|\varphi_i(x)|)^m}
\bigl|F_{\alpha,\lambda}^{\varphi_i,\psi_i}(x)\bigr|
\le CC_m<\infty,
\]
which contradicts~\cite[(2.3)]{AJK}. This proves the lemma.
\end{proof}

\section{Amendment to Theorem 2.11}
Theorem~2.11 is false as stated. Nevertheless, under the additional hypothesis that $\varphi$ is bounded, the proof is correct and, under this additional assumption, is simplified. Accordingly, the statement is replaced here by its bounded-symbol version:

\medskip
\noindent\emph{Let $\varphi:\R^d\to\R^d$ be smooth and bounded.
Then $C_{\psi,\varphi}$ acts on $\Sch(\R^d)$ for every
$\psi\in\Sch(\R^d)$ if and only if
$\varphi_j\in\mathcal O_M(\R^d)$ for all
$1\le j\le d$.}
\medskip

The restriction to bounded symbols is not intended to be optimal. For
extensions and related statements we refer to \cite{LiCorrections}.

\section{Amendment to Lemma 3.1}
The proof of Lemma~3.1 contains the same unjustified assumption as the original proof of Lemma~2.1, namely that the points witnessing the unboundedness in condition~(3.3) of \cite{AJK} can be chosen arbitrarily far from the origin. Unlike Lemma~2.1, however, Lemma~3.1 is false as stated. Its conclusion nevertheless remains valid if condition~(3.3) is replaced by the following stronger condition: there is $\alpha\in\N_0^d$ such that
\begin{equation}
\limsup_{R\to\infty}\ \sup_{i\in I}\sup_{|x|>R}
\frac{|\varphi_i^{(\alpha)}(x)|}
     {(1+|\varphi_i(x)|)^q}
=\infty
\tag{3.3'}\label{eq:modified33}
\end{equation}
for all $q>0$.

Indeed, following the original proof up to equation~(3.4) in \cite{AJK}, we choose $\alpha\in\N_0^d$ with minimal $|\alpha|$ such that \eqref{eq:modified33} holds for every $q>0$. Again, $|\alpha|\geq 1$, and for every $\beta\in\N_0^d$ with $\beta\leq \alpha$, $\beta\neq 0$, there is $q'(\beta)>0$ such that
\begin{equation}
\limsup_{R\to\infty}\ \sup_{i\in I}\sup_{|x|>R}
\frac{|\varphi_i^{(\alpha-\beta)}(x)|}
     {(1+|\varphi_i(x)|)^{q'(\beta)}}
<\infty.
\tag{3.5'}\label{eq:modified35}
\end{equation}
The remainder of the proof proceeds as in \cite{AJK} with the points $x_j$ chosen sufficiently far from the origin so that the bounds provided by \eqref{eq:modified35} apply. Thus, the global estimates~(3.5) in \cite{AJK} are replaced by the corresponding estimates valid outside a sufficiently large ball.

We emphasize that conditions~(3.1) and~(3.3) of \cite[Lemma~3.1]{AJK} do not, in general, imply condition~\eqref{eq:modified33}. Consequently, the amendment does not restore the original statement or its subsequent applications without further verification.

However, we explicitly point out that, if $I=\{1\}$, then condition~(3.3) of \cite{AJK} and condition~\eqref{eq:modified33} are equivalent. In particular, Theorem~3.3 of \cite{AJK} remains valid.

\section{Amendments to Sections 4 and 5}
The amendment to Lemma~3.1 affects several results in Sections~4 and~5 of \cite{AJK}. In the following, we describe the necessary modifications, provide corrected proofs where needed, and withdraw those assertions that cannot be maintained.
\begin{enumerate}[1.]
\item \textbf{Theorem 4.5.} Due to the amendment to Lemma~3.1 and the fact that $\psi^{n,\varphi}$ need not be of small decay with respect to $\varphi_n$, even when $\psi$ is of small decay with respect to $\varphi$, condition~(ii)(b) of Theorem~4.5 does not, in general, follow from condition~(i). Theorem~4.5 must therefore be replaced by the following statement, which is precisely the assertion of Remark~4.6 of \cite{AJK}. We include a proof for the reader's convenience.

\emph{Let $\varphi\colon\R^d\rightarrow\R^d$ be a smooth mapping and let $\psi\in C^\infty(\R^d)$. Then, for the following conditions, (ii) implies (i), and (i) implies (ii)(a).
\begin{enumerate}[(i)]
\item $C_{\psi,\varphi}$ acts on $\mathscr{S}(\R^d)$ and is power bounded.
\item We have the following two conditions.
\begin{enumerate}[(a)]
\item For all $p>0$ and $\alpha\in\N_0^d$ there exists $q>0$ such that
\[\sup_{n\in\N}\sup_{x\in\R^d}\frac{(1+|x|)^p}{(1+|\varphi_n(x)|)^q}|(\psi^{n,\varphi})^{(\alpha)}(x)|<\infty.\]
\item For all $\alpha\in\N_0^d$ there exists $q>0$ such that
\[\sup_{n\in\N}\sup_{x\in\R^d}\frac{1}{(1+|\varphi_n(x)|)^q}|\varphi_n^{(\alpha)}(x)|<\infty.\]
\end{enumerate} 
\end{enumerate}
} 

\emph{Proof.}
Suppose that (ii)(a) and (ii)(b) hold. By Corollary~2.4 of \cite{AJK}, $C_{\psi,\varphi}$ acts on $\mathscr{S}(\R^d)$. Moreover,~(2.1) and~(2.2) of~\cite{AJK} imply that $\{C_{\psi,\varphi}^n f:n\in\N\}$ is bounded in $\mathscr{S}(\R^d)$ for every $f\in\mathscr{S}(\R^d)$. Thus, the Uniform Boundedness Principle yields the power boundedness of $C_{\psi,\varphi}$.

Conversely, suppose that~(i) holds. Since $F^{\varphi_n,\psi^{n,\varphi}}_{\alpha,0}
=(\psi^{n,\varphi})^{(\alpha)}$ by equation~(2.1) of~\cite{AJK}, Lemma~2.1 applied to the family $(\psi^{n,\varphi})_{n\in\mathbb{N}}$ and to the smooth symbols $(\varphi_n)_{n\in\N}$ yields~(ii)(a).
\hfill$\square$\\

\item \textbf{Theorem 4.12.}
For the same reasons as in Theorem~4.5, condition~(ii)(b) of Theorem~4.12 does not, in general, follow from condition~(i). Thus, Theorem~4.12 must be replaced by the following statement, which is precisely the assertion of Remark~4.13 of \cite{AJK}. Again, we include a proof for the reader's convenience.

\emph{Let $\varphi\colon\R^d\to\R^d$ be a smooth mapping and let $\psi\in C^\infty(\R^d)$. Then, for the following conditions, (ii) implies (i) and (i) implies (ii)(a).
\begin{enumerate}[(i)]
\item $C_{\psi,\varphi}$ acts on $\mathscr{S}(\R^d)$ and is topologizable (m-topologizable).
\item We have the following to conditions.
\begin{enumerate}[(a)]
\item For all $p>0$ and $\alpha\in\N_0^d$ there exists $q>0$ such that, for every $n\in\N$,
\[
\sup_{x\in\R^d}
\frac{(1+|x|)^p}{(1+|\varphi_n(x)|)^q}
|(\psi^{n,\varphi})^{(\alpha)}(x)|
<\infty
\quad
(\leq M^n\text{ for some }M>0).
\]
\item For all $\alpha\in\N_0^d$ there exists $q>0$ such that, for every $n\in\N$,
\[
\sup_{x\in\R^d}
\frac{|\varphi_n^{(\alpha)}(x)|}
     {(1+|\varphi_n(x)|)^q}
<\infty
\quad
(\leq M^n\text{ for some }M>0).
\]
\end{enumerate}
\end{enumerate}
}

\emph{Proof.}
Suppose that (ii)(a) and (ii)(b) hold. By Corollary~2.4 of \cite{AJK}, $C_{\psi,\varphi}$ acts on $\mathscr{S}(\R^d)$. Moreover, equations~(2.1) and~(2.2) of \cite{AJK} imply that, for every $p>0$ and $\alpha\in\N_0^d$, there exists $q>0$ such that
\[
\forall n\in\N\ \exists C_n>0\ \forall f\in\mathscr{S}(\R^d):
\quad
\|C_{\psi,\varphi}^n f\|_{p,\alpha}
\leq C_n\max_{|\lambda|\leq|\alpha|}
\|f\|_{q,\lambda},
\]
where the constants $C_n$ can be chosen of the form $M^n$ for some $M>0$ in the m-topologizable case. Hence, $C_{\psi,\varphi}$ is topologizable (m-topologizable).

Conversely, suppose that (i) holds. By Proposition~4.10 of \cite{AJK}, for every $p>0$ and $\alpha\in\N_0^d$ there exists $q>0$ such that
\[
\sup_{x\in\R^d}
\frac{(1+|x|)^p}{(1+|\varphi_n(x)|)^q}
|F^{\varphi_n,\psi^{n,\varphi}}_{\alpha,0}(x)|
<\infty
\quad
(\leq M^n\text{ for some }M>0)
\]
for every $n\in\N$. Since $F^{\varphi_n,\psi^{n,\varphi}}_{\alpha,0}
=(\psi^{n,\varphi})^{(\alpha)}$
by equation~(2.1) of~\cite{AJK}, condition~(ii)(a) follows.
\hfill$\square$\\


\item \textbf{Proposition 5.4.}
Since condition~(ii)(b) of Theorem~4.5 is not necessary for the power boundedness of $C_{\psi,\varphi}$ on $\mathscr{S}(\R^d)$, the proof of Proposition~5.4 in \cite{AJK} is invalid for $|a|>1$. We therefore restrict Proposition~5.4 to the following statement, whose proof is already contained in \cite{AJK}.

\emph{Let $\varphi\colon\R\to\R$, $x\mapsto ax+b$, with $b\in\R$ and $0<|a|<1$, and let $\psi$ be a nonzero polynomial. Then $C_{\psi,\varphi}$ is not power bounded.}

We note that the conclusion also holds for $|a|>1$ whenever $\psi\left(\frac{b}{1-a}\right)\neq 0$, as observed in \cite{LiCorrections}.\\

\item \textbf{Proposition 5.5.} Proposition~5.5 is false as stated and must be withdrawn. We refer to \cite{LiCorrections} for related results.\\


\item \textbf{Theorem 5.8.} In condition~(ii), remove the alternative concerning smooth weights $\psi$ of small decay with respect to $\varphi$. Thus, condition~(ii) must be replaced by the requirement that $C_{\psi,\varphi}$ is power bounded on $\mathscr{S}(\R^d)$ for every $\psi\in\mathcal O_M(\R^d)$.\\

\item \textbf{Theorem 5.10.}
In conditions~(f), (g), and~(h), remove the alternative concerning smooth weights $\psi$ of small decay with respect to $\varphi$. Thus, these conditions must be replaced by the following:

\begin{enumerate}
\item[(f)] $C_{\psi,\varphi}$ is power bounded on $\Sch(\R)$ for every $\psi\in\mathcal O_M(\R)$.
\item[(g)] $(C_{\psi,\varphi}^n)_{n\in\N}$ converges to $0$ in $\mathcal L_b(\Sch(\R))$ for every $\psi\in\mathcal O_M(\R)$.
\item[(h)] $C_{\psi,\varphi}$ is (uniformly) mean ergodic on $\Sch(\R)$ for every $\psi\in\mathcal O_M(\R)$.
\end{enumerate}

Moreover, the original proof of the implication~(e)$\Rightarrow$(b) relies on Proposition~5.5 and is therefore invalid. To complete the proof of the corrected theorem, we establish the implication~(e)$\Rightarrow$(a) directly, using Ces\`aro means, and recall that the equivalence of conditions~(a)--(d) has been established in the literature; see~\cite{AJK} for the references. We write
\[
 A_nf=\frac1n\sum_{m=1}^n C_\varphi^m f.
\]
Assume (e), i.e.~$\{A_n\colon n\in\N\}$ is bounded in $\mathcal{L}_b(\mathscr{S}(\R))$, and suppose, to derive a contradiction, that~(a) does not hold, i.e.~that $\varphi$ has a fixed point. Let $x_0$ be its largest fixed point. Since $\varphi(x)-x$ has no zeros on $(x_0,\infty)$, it has constant sign there.

\emph{Case 1: $\varphi(x)>x$ for every $x>x_0$.} For every such $x$, the sequence $\varphi^n(x)$ is strictly increasing. It cannot be bounded, since its limit would be a fixed point greater than $x_0$. Hence $\varphi_n(x)\to+\infty$. For every $f\in\Sch(\R)$ it follows that
\begin{equation}\label{eq:cesaromeans}
 A_nf(x)\to 0\quad(x>x_0),\qquad
 A_nf(x_0)=f(x_0).
\end{equation}
By (e), the sequence $(A_nf)_{n\in\mathbb{N}}$ is bounded in $\mathscr{S}(\R)$. Since $\mathscr{S}(\R)$ is a Montel space, $(A_nf)_{n\in\mathbb{N}}$ has a convergent subsequence. Its limit is a continuous function, contradicting the pointwise limits in~\eqref{eq:cesaromeans} for a choice of $f$ with $f(x_0)\ne0$.

\emph{Case 2: $\varphi(x)<x$ for every $x>x_0$.}
Then $\varphi$ has negative leading coefficient. Set
$\theta=\varphi_2$. We distinguish the parity of $\deg\varphi$.

If $\deg\varphi$ is odd, then $\theta$ has odd degree and positive
leading coefficient. Let $x_1$ be its largest fixed point, which
exists since $x_0$ is a fixed point of $\theta$. We have
$\theta(x)>x$ for $x>x_1$. Thus $\theta_n(x)\to+\infty$
for all $x>x_1$, by the monotonicity argument in Case~1.
Since $\varphi(t)\to-\infty$ as $t\to+\infty$, the odd
iterates $\varphi_{2n+1}(x)$ also tend to infinity in modulus.
Consequently $A_nf(x)\to0$ for $x>x_1$ and every
$f\in\Sch(\R)$.

If $\deg\varphi$ is even, then $\theta$ has even degree
and negative leading coefficient. Let $x_1$ now be the
\emph{smallest} fixed point of $\theta$. Since
$\theta(x)-x\to-\infty$ as $x\to-\infty$ and there
are no fixed points of $\theta$ below $x_1$, we have
$\theta(x)<x$ for $x<x_1$. The sequence
$\theta_n(x)$ is strictly decreasing and must tend to
$-\infty$ for every $x<x_1$; otherwise its finite limit
would be a fixed point below $x_1$. Again the odd iterates
$\varphi_{2n+1}(x)=\varphi(\theta_n(x))$ tend to infinity
in modulus. Hence $A_nf(x)\to0$ for $x<x_1$.

In either subcase, $\theta(x_1)=x_1$, so the orbit of $x_1$
under $\varphi$ has period at most two and
\[
A_nf(x_1)\to
\frac{f(x_1)+f(\varphi(x_1))}{2}.
\]
Choose $f\in\Sch(\R)$ such that
$f(x_1)=f(\varphi(x_1))=1$.
The pointwise limit of $(A_nf)_{n\in\N}$ is therefore
discontinuous at $x_1$, since it vanishes on an adjacent
half-line. On the other hand, by (e), $(A_nf)_{n\in\N}$
is bounded in the Montel space $\Sch(\R)$ and hence
has a convergent subsequence whose limit is continuous.
This contradiction proves (e)$\Rightarrow$(a).

\end{enumerate}

\bibliographystyle{plain}
\bibliography{references}

\end{document}